%% file: agt-4-38.tex
\documentclass{gtart_h}

\input agtout

\lognumber{38}
\volumenumber{4}
\volumeyear{2004}
\papernumber{38}
\pagenumbers{861}{892}
\received{17 September 2003} 
\revised{30 July 2004}
\accepted{13 September 2004}
\published{9 October 2004}

\usepackage{amsmath,amssymb}
\let\Bbb\mathbb

\theoremstyle{plain}
\newtheorem{thm}{Theorem}[section]
\newtheorem{lem}{Lemma}[section]

\newtheorem{cor}{Corollary}[section]
\newtheorem{prop}{Proposition}[section]

\theoremstyle{definition}

\newtheorem{prob}{Problem}[section]

\newtheorem{question}{Question}[section]

\theoremstyle{remark}
\newtheorem{remark}{Remark}[section]

\numberwithin{equation}{section}

\DeclareMathOperator{\cat0dim}{CAT(0)-dim}
\DeclareMathOperator{\geomdim}{geom \, dim}
\DeclareMathOperator{\dimss}{CAT(0)-dim_{ss}}
\DeclareMathOperator{\cd}{cd}

\DeclareMathOperator{\CAT}{CAT}
\DeclareMathOperator{\Isom}{Isom}

\newcommand{\field}[1]{\mathbb{#1}}

\newcommand{\R}{\field{R}}

\newcommand{\Z}{\field{Z}}

\newenvironment{proclaim}[1]
{\par\medskip\noindent{\bf #1}\qua \sl}
{\rm\par\medskip}

\title[Parabolic isometry and CAT(0) dimension]
{Parabolic isometries of CAT(0) spaces\\and CAT(0) dimensions}
\author{Koji Fujiwara\\Takashi Shioya\\Saeko Yamagata}
\shortauthors{Fujiwara, Shioya and Yamagata}
\gtemail{\mailto{fujiwara@math.tohoku.ac.jp}, 
\mailto{shioya@math.tohoku.ac.jp}, \mailto{sa1m28@math.tohoku.ac.jp}}
\asciiemail{fujiwara@math.tohoku.ac.jp, shioya@math.tohoku.ac.jp, sa1m28@math.tohoku.ac.jp}
\address{Mathematics Institute, Tohoku University, Sendai 980-8578, Japan}



\begin{abstract}
We study discrete groups from the view point of
a dimension gap in connection to CAT(0) geometry.
Developing studies by Brady-Crisp and Bridson,
we show that there exist finitely presented groups of geometric
dimension $2$
which do not act properly on any proper $\CAT(0)$ spaces
of dimension $2$ by isometries, although such actions exist
on $\CAT(0)$ spaces of dimension $3$.

Another example is the fundamental group, $G$, of a complete, non-compact,
complex hyperbolic manifold $M$ with finite volume, of complex-dimension
$n \ge 2$. The group $G$ is acting on the universal cover of $M$, which
is isometric to ${\bf H}^n_{{\Bbb C}}$. It is a $\CAT(-1)$ space of
dimension $2n$.
The geometric dimension of $G$ is $2n-1$.
We show that $G$ does not act on any proper $\CAT(0)$ space of dimension
$2n-1$
properly by isometries.

We also discuss the fundamental groups of a torus bundle over a circle,
and solvable Baumslag-Solitar groups.
\end{abstract}

\asciiabstract{%
We study discrete groups from the view point of a dimension gap in
connection to CAT(0) geometry.  Developing studies by Brady-Crisp and
Bridson, we show that there exist finitely presented groups of
geometric dimension 2 which do not act properly on any proper
CAT(0) spaces of dimension 2 by isometries, although such actions
exist on CAT(0) spaces of dimension 3.

Another example is the fundamental group, G, of a complete,
non-compact, complex hyperbolic manifold M with finite volume, of
complex-dimension n > 1.  The group G is acting on the universal cover
of M, which is isometric to H^n_C. It is a CAT(-1) space of dimension
2n.  The geometric dimension of G is 2n-1.  We show that G does not
act on any proper CAT(0) space of dimension 2n-1 properly by
isometries.

We also discuss the fundamental groups of a torus bundle over a circle,
and solvable Baumslag-Solitar groups.}

\begin{document}
\maketitle
\newtheorem{theorem}{Theorem}[section]

\newcommand{\G}{\Gamma}

\section{Introduction and statement of results}

As a generalization of simply connected Riemannian manifolds
of non-positive sectional curvature, the notion
of ``CAT(0)" was introduced to geodesic spaces by M.Gromov.
Many properties of simply connected
Riemannian manifolds of non-positive sectional curvature
remain true for CAT(0) spaces. We refer to \cite{BB}
and also \cite{BH} for a thorough account of the subject.
For example, if a space is CAT(0) then it is contractible,
which is a Cartan-Hadamard theorem.
This gives a powerful tool to show some spaces are contractible
(for example, see \cite{CD}).
Another example is the classification of isometries of a complete
CAT(0) space $X$: elliptic, hyperbolic, or parabolic.
It seems harder to understand parabolic isometries
than the other two, which give natural objects to look at:
the set of fixed points for an elliptic
isometry and the set of axes for a hyperbolic isometry.
Moreover the union of the axes is isometric to the product
of a convex set in $X$ and a real line.
``Flat torus theorem" (see Theorem \ref{flattorus}) is an important
consequence
of the existence of axes.
An isometry is called {\it semi-simple} if it is either
elliptic or hyperbolic.

A metric space, $X$, is called {\it proper} if for any $x \in X$ and any
$r >0$,
the closed ball in $X$ of radius $r$ centered at $x$ is compact.
If $X$ is proper, then it is locally compact.
If $X$ is a complete, locally compact, geodesic space,
then it is proper (cf.\ \cite{BuBuI} Prop 2.5.22).
Therefore, a complete CAT(0) space is proper
if and only if it is locally compact.

If $X$ is a proper CAT(0) space, then a parabolic isometry
has at least one fixed point in the ideal boundary of $X$
(see Proposition \ref{parabolicfix}).

Suppose a discrete group $G$ is acting on a topological space, $X$,
by homeomorphisms.
To give definitions of ``proper actions", let's consider
the following conditions.
\begin{enumerate}
\item
For any $x \in X$, there exists a neighborhood $U \subset X$ of $x$
such that $\{ g \,|\, U \cap gU \not = \emptyset \}$ is finite.
\item
For any $x, y \in X$, there exist neighborhoods $U$ of $x$ and
$V$ of $y$ such that
$\{ g \,|\, V \cap gU \not = \emptyset \}$ is finite.
\item
For any compact subset $K \subset X$,
$\{ g \,|\, K \cap gK \not = \emptyset \}$ is finite.
\item
The map $G \times X \to X \times X$, $(g,x) \mapsto (gx,x)$
is a proper map.
\end{enumerate}

Recall that a continuous map, $f:X\to Y$, is
called {\it proper} if it is closed and for each $y \in Y$ ,
$f^{-1}(y)$ is compact.
If $X$ and $Y$ are Hausdorff and $Y$ is locally compact,
then $f$ is proper if and only if for each compact set
$K \subset Y$, the preimage $f^{-1}(K)$ is compact (\cite{Di} Ch1 ).

Clearly (2) implies (1). It is easy to see
that (2) implies (3).
If $X$ is locally compact, then (3) implies (2).
If $X$ is Hausdorff, then (2) and (4) are equivalent
(\cite{Di} 3.22 Cor).
If $X$ is a metric space and the action is
by isometries, then (1) and (2)
are equivalent (\cite{DV} Thm1),
so that (1),(2) and (4) are equivalent, and (3) follows
from one of them.

Following \cite{BH}, we say that an isometric action of a discrete
group, $G$, on
a metric space, $X$, is {\it proper} if the
condition (1) is satisfied, in other words, for any $x$, there exists
$r>0$ such that
$\{ g\,|\,g B(x,r) \cap B(x,r) \not = \emptyset\}$ is finite, where
$B(x,r) =\{ y \in X\,|\, d(x,y) \le r \}$.
In the literature, proper actions are often called
properly discontinuous.
Most of the results in this paper are stated
for proper actions, which is our main theme,
but sometimes
the weaker condition (3) is enough.
Since we can not find a standard name,
we decide to call an action with the condition
(3) a {\it K-proper} action. As we said, if
the space is locally compact, then those two
properness are equivalent.

%
%

In this paper we study proper, isometric actions of a discrete group $G$ on
CAT(0) spaces $X$.
We are interested in the minimal dimension of such $X$, which
is called a {\it CAT(0) dimension} of $G$, denoted CAT$(0)$-dim $G$.
If such actions do not exist, then CAT(0)-dim $=\infty$.

By the dimension of a topological (or
metric) space, we mean
the covering dimension (see subsection \ref{covdim} for the precise
definition).

If $G$ is torsion-free, there is an obvious lower bound of $\cat0dim G$.
Recall that the {\it geometric dimension} $\le \infty$ of $G$,
denoted $\geomdim G$,
is the minimal dimension of a $K(G,1)$-complex.
If $G$ is torsion-free, then a proper action of $G$ on
a CAT(0) space $X$
is free, so that $X/G$ is a $K(G,1)$-space,
because $X \to X/G$ is a covering map.
It follows that $\geomdim G \le \cat0dim G$.

We are interested in the following problem.
\begin{prob}\label{problem}
Find a torsion-free group $G$ such that
$$\geomdim G < \cat0dim G < \infty.$$
\end{prob}

It seems an example of such $G$ is not
known, but a few potential examples of groups are studied
by Brady-Crisp \cite{BC} and Bridson \cite{Br}
under an extra assumption on the actions.
The following conditions on
actions of $G$ on a CAT(0) space $X$
are natural.

\begin{itemize}
\item[(cc)]The action is co-compact, i.e., $X/G$ is compact.

\item[(ss)]The action is by semi-simple isometries.

\item[(p)]$X$ is proper.
\end{itemize}

A group which acts properly on some CAT(0) space
with the condition (cc) is called a {\it CAT(0) group}.
The condition (cc) implies (ss).

If we choose one (or more) of the above conditions and only consider
proper, isometric actions of $G$ on CAT(0) spaces
which satisfy the condition, we obtain another definition
of a CAT(0) dimension of $G$, which is clearly not smaller than
the original one.
The condition (ss) is the one Brady-Crisp and
Bridson imposed.
Bridson denotes the CAT(0) dimension in this sense by
$\dimss$.
As for this dimension, in other words,
if we consider only actions of $G$ on CAT(0) spaces
by semi-simple isometries, it is easy
to find $G$ such that $\geomdim$ is finite and
$\dimss$ is $\infty$ (for example, see Proposition \ref{bs}).
The important part of their work is that
they found groups with a finite gap, which is $1$.
Their examples are CAT(0) groups.

In this paper we study a group action on a CAT(0) space
such that the action has a fixed point in the ideal boundary.
It turns out that under certain circumstances, the condition
(p) implies (ss).
As an application, we replace the assumption
(ss) by (p) in the results of Brady-Crisp and Bridson.

To analyze an action with a common fixed point in the
ideal boundary, the following result, (2), on the
dimension of horospheres is essential.
Let $X(\infty)$ be the ideal boundary of a CAT(0) space, $X$, and
$\bar{X}=X \cup X(\infty)$ with the cone topology (see section \ref{def} for
definitions).
A metric sphere is defined by
$S(x,r)= \{ y \in X\,|\, d(x,y) =r\}$.

\begin{proclaim}{Theorem \ref{bestvina}.}
Let $X$ be a proper $\CAT(0)$ space of dimension
$n < \infty$. Then
\begin{enumerate}
\item
A metric sphere $S(x,r) \subset X$
has dimension at most $n-1$.
\item
A horosphere centered at a point $p \in X(\infty)$
has dimension at most $n-1$.
\item
$\bar{X}$ is an AR of dimension $n$.
The ideal boundary $X(\infty)$ is a $Z$-set in $\bar{X}$,
and the dimension is at most $n-1$.
$\bar{X}$ is homotopy equivalent to $X$, hence
contractible.
\end{enumerate}
\end{proclaim}

Let $\cd G$ denote the cohomological dimension of $G$
(the definition is in subsection \ref{cohodim}).
We show the next result using Theorem \ref{bestvina}(2).
Although the dimension of a horosphere is at most
$n-1$ by (2), we need some extra work
to show $\cd G \le n-1$ because
the horosphere is not contractible in general.

\begin{proclaim}{Proposition \ref{cdG}.}
Let $X$ be a proper $\CAT(0)$ space of dimension $n < \infty$.
Suppose a group
$G$ is acting on $X$ properly by isometries.
Assume $G$ fixes a point $p \in X(\infty)$ and leaves each horosphere $H_t$
centered at $p$ invariant.
If $G$ has a finite $K(G,1)$-complex then
$\cd G \le n-1$.
In particular, if $n=2$ then $G$ is free.
\end{proclaim}

It has a useful consequence.

\begin{proclaim}{Proposition \ref{ss}.}
Let $X$ be a proper $\CAT(0)$ space.
Suppose a free abelian group of rank $n$, ${\Bbb Z}^n$,
is acting properly on $X$ by isometries.
Then $\dim X \ge n$.
If $\dim X =n$
then the action is by semi-simple isometries.
\end{proclaim}

Combining this proposition and the study by Brady-Crisp \cite{BC} or Bridson
\cite{Br} we obtain the following result, which
is an immediate consequence of
Theorems \ref{BC} and \ref{BC2}, or Theorems \ref{bridson} and
\ref{bridson2}.

\begin{thm}\label{cat0}
There exists a finitely presented group $G$
of geometric dimension $2$ such that

\begin{enumerate}
\item
$G$ does not act properly on any proper $\CAT(0)$ space of dimension
$2$ by isometries.

\item
$G$ acts properly and co-compactly on some proper $\CAT(0)$ space of
dimension $3$ by isometries.
In particular, $G$ is a $\CAT(0)$ group.
\end{enumerate}
\end{thm}

Also we give a new class of examples which also
have gaps between the geometric dimension and the
CAT(0) dimension with the condition (p).
Those groups are not CAT(0) groups.
See Cor \ref{heisenberg} for the class of examples
of various dimensions and a proof.
Let ${\bf H}_{\Bbb C}^n$ denote
the complex hyperbolic space of complex-dimension
$n$. It is a complete, simply connected Riemannian
manifold of (real-) dimension $2n$ with the
sectional curvature pinched by $-1$ and $-1/4$, hence
a CAT$(-1)$ space.

\begin{thm}
Let $G=\langle a,b,c\,|\,[a,b]=c,[a,c]=1,[b,c]=1 \rangle$.
Then
\begin{enumerate}
\item
$\geomdim G=3$.
\item
$G$ does not act freely, hence not properly,
on any $\CAT(0)$ space
by semi-simple isometries,
so that $G$ is not a $\CAT(0)$ group.
\item
$G$ does not act properly on any proper $\CAT(0)$ space
of dimension $3$ by
isometries.

\item
$G$ acts on ${\bf H}_{\Bbb C}^2$, which is a proper $\CAT(-1)$ space
of dimension $4$, properly by
isometries.
\end{enumerate}
\end{thm}

We also show the following.

\begin{proclaim}{Corollary \ref{complex}.}
Let $G$ be the fundamental group of a non-compact, complete,
complex hyperbolic manifold $M$
of complex-dimension $n \ge 2$. If the volume of $M$ is finite, then

\begin{enumerate}
\item
$\cd G=\geomdim G=2n-1$.

\item
$G$ does not act properly on
any proper $(2n-1)$-dimensional $\CAT(0)$ space
by isometries.
\item
$G$ acts properly on ${\bf H}_{\Bbb C}^{n}$, which is a proper
$2n$-dimensional $\CAT(-1)$ space, by isometries.
\item
$G$ does not act on any $\CAT(0)$ space
freely, hence not properly, by semi-simple isometries,
so that $G$ is not a $\CAT(0)$ group.

\end{enumerate}
\end{proclaim}

We ask several questions in the paper, including the
following two (see discussions later).

Let $BS(1,m)=\langle a,b\,|\,a b a^{-1}=b^m \rangle.$
$BS(1,m)$ is called a solvable Baumslag-Solitar group.

\medskip
{\bf Question \ref{bs.question}}\qua
Does $BS(1,m)$ act properly
on a CAT(0) space $X$ of dimension $2$?\medskip

Let
$S=\langle a,b,c\,|\,ab=ba,cac^{-1}=a^2b,cbc^{-1}=ab \rangle.$
$S$ is the fundamental group of a closed
$3$-manifold which is a torus
bundle over a circle.

\noindent
{\bf Question \ref{bundle}}\qua
Does $S$ act properly on some $3$-dimensional CAT(0) space
by isometries?\medskip

\noindent
{\bf Acknowledgement}\qua
We would like to thank A.Casson, H.Izeki, Y.Kamishima, M.Kanai,
K.Nagano, K.Ono, and M.Sakuma for useful comments and information.
We thank the referee for comments which improved
the presentation of the paper.

The first author owes Wasatch Topology Conference
in Park City and LMS Durham Symposium in the Summer 2003,
where he had valuable discussions with several
participants; in particular, J.Crisp,
N.Monod, E.Swenson and K.Whyte.
He is grateful to M.Bestvina for his constant interest
and suggestions.

\section{A CAT(0) space and its ideal boundary}\label{def}
In this section we review basic definitions and facts
on CAT(0) spaces. See \cite{BH} for details.

\subsection{CAT(0) spaces}
Let $X$ be a geodesic space and $\Delta(x,y,z)$
a geodesic triangle in $X$, which is a union
of three geodesics.
A {\it comparison triangle} for $\Delta$ is a triangle
$\bar\Delta(\bar x,\bar y,\bar z)$ in ${\bf E}^2$
with the same side lengths as $\Delta$.
Let $p$ be a point in the side of $\Delta$ between, say,
$x,y$, which is often denoted by $[x,y]$.
A {\it comparison point} in $\bar \Delta$ is a point $\bar p \in [\bar
x,\bar y]$
with $d(x,p)=d_{{\bf E}^2}(\bar x ,\bar p)$.
$\bar \Delta$ {\it satisfies the CAT(0) inequality}
if for any $p,q \in \Delta$ and their comparison points
$\bar p,\bar q \in \bar \Delta$, $d(p,q)=d_{{\bf E}^2}(\bar p,\bar q)$.
$X$ is a {\it CAT(0) space} if all geodesic triangles in $X$
satisfy CAT(0) inequality.

Similarly, one defines a notion of {\it CAT$(1)$} and {\it CAT$(-1)$} spaces
by comparing every geodesic triangle in $X$
with its comparison
triangle in the standard $2$-sphere ${\bf S}^2$ and the
real $2$-dimensional hyperbolic space ${\bf H}^2$,
respectively. In the case of CAT(1) we only consider geodesic triangles
of total
perimeter length less than $2 \pi$.

\subsection{Ideal boundaries}
Two geodesics $\gamma(t),\gamma'(t)$ in a complete
CAT(0) space $X$ are {\it asymptotic}
if there exists a constant $C$ such that for all $t \ge 0$,
$d(\gamma(t),\gamma'(t)) \le C$.
This defines an equivalence relation, $\sim$.
In this paper all geodesics have unit speed.
The set of {\it points at infinity}, $X(\infty)$,
or the {\it ideal boundary} of $X$
is the set of equivalence classes of geodesics in $X$.
The equivalence class of a geodesic $\gamma(t)$ is denoted by
$\gamma(\infty)$. The equivalence class of a geodesic $\gamma(-t)$
is denoted by $\gamma(-\infty)$.
It is a theorem that for any geodesic $\gamma(t)$ and a point $x \in X$,
there exists a unique geodesic $\gamma'(t)$ such that $\gamma \sim \gamma'$
and $\gamma'(0)=x$.

\subsection{The cone topology}
Let $X$ be a complete CAT(0) space and $\bar{X}=X \cup X(\infty)$.
One standard topology on $\bar{X}$ is called
the {\it cone topology}. The induced topology on $X$
is the original one.
We will give a basis for the cone topology.
A typical basis for $x \in X$ is
$B(x,r) \subset X$.
To give a typical basis for a point $p \in X(\infty)$, fix
$x \in X$ and $r>0$. Let $c$ be the geodesic such that
$c(0)=x$ and $c(\infty)=p$.
We extend the nearest point projection
${\rm pr}:X \to B(x,r)$ to $\bar{X}$.
For any point $q \in X(\infty)$ let $\alpha_q$
be the geodesic with $\alpha_q(0)=x$ and $\alpha_q(\infty)=q$.
Define ${\rm pr}(q) \in B(x,r)$ to be $\alpha_q(r)$.
For $\epsilon >0$ let
$$U(x,r,\epsilon;p)=\{x \in \bar{X}\,|\,d(c(r),pr(x)) < \epsilon \}.$$
There are two types of explicit neighborhood bases
for the cone topology
on $\bar{X}$: all $B(x,r)$ for $x \in X$ and $r>0$, and
all $U(x,r,\epsilon;p)$ for $x \in X, r>0, \epsilon>0, p \in X(\infty)$.
We remark that not only $B(x,r)$ but also $U(x,r,\epsilon,p) \cap X$
is convex in $X$.

Let $f$ be an isometry of $X$. Then since $f$
acts on the set of geodesics on $X$ leaving the
equivalence relation $\sim$ invariant, it gives
an action on $X(\infty)$. It is known this action
of $f$ on $\bar{X}$ is a homeomorphism.

If $X$ is a proper CAT(0) space, then $\bar{X}$
is compact with respect to the cone topology.
Also it is metrizable since it is a second-countable,
normal space (Urysohn's Metrization Theorem).

There exists a (pseudo)metric on $X(\infty)$
which is called the {\it Tits metric}. The topology induced
by the Tits metric is generally stronger than the cone
topology. $\bar{X}$ may not be compact with respect to
the Tits topology. This is one reason why we
only use the cone topology in this paper.

\subsection{Busemann functions and horospheres}

Let $X$ be a complete CAT(0) space and let
$c:[0,\infty] \to X$ be a geodesic ray.
The function $b_c:X \to {\Bbb R}$
defined by the following is called the {\it Busemann function}
associated to $c$:
$$b_c(x)=\lim_{t \to \infty}(d(x,c(t))-t).$$
$b_c$ is a convex function such that
for all $x,y \in X$, $|b_c(x)-b_c(y)| \le d(x,y)$.
It is known (see 8.20 Cor in \cite{BH}) that
if $c \sim c'$ then there exists a constant $C$
such that for all $x \in X$, $b_c(x)-b_{c'}(x)=C$.
Because of this, level sets of $b_c$,
$\{x \in X\,|\,b_c(x)=t\} \subset X, t \in {\Bbb R}$, depend
only on the asymptotic class of $c$. They are called
{\it horospheres} (centered) at $c(\infty)$.
By construction of horospheres an isometry $f$
of $X$ which fixes $c(\infty)$ acts on the set
of horospheres at $c(\infty)$.
A {\it horoball} is defined by
$\{ x \in X\,|\,b_c(x) \le t \}$. Horoballs are convex subsets
in $X$.

\section{Parabolic isometries}
\subsection{Classification of isometries}
Let $f$ be an isometry of a complete CAT(0) space $(X,d)$.
The {\it displacement function}
$d_f:X \to {\Bbb R}$ is defined
by
$$d_f(x)=d(x,f(x)).$$
This is a convex, $f$-invariant function.
There is the following classification of isometries $f$ (cf.\ \cite{BH}).
\begin{itemize}
\item
If $f$ fixes a point in $X$ it is called {\it elliptic}.
\item
If there is a bi-infinite geodesic $\gamma$ in $X$ which
is invariant by $f$ and $f$ acts on it by non-trivial
translation, then $f$ is called {\it hyperbolic}.
$\gamma$ is called {\it axis}.
An axis may not be unique, but they are parallel to each other,
i.e., there is a convex subspace between them
which is isometric the product of an interval and ${\bf R}$,
called a {\it strip}.
\item
Otherwise $f$ is called {\it parabolic}.
\end{itemize}

If $f$ is elliptic or hyperbolic,
it is called {\it semi-simple}.
It is known that $f$ is hyperbolic if and only if $\inf d_f >0$
and the infimum is attained.
A point $x \in X$ is on an axis if and only if $d_f$ attains
its infimum at $x$.
$f$ is parabolic if and only if $d_f$ does not attain its infimum.
The infimum may be positive and sometimes $f$ is called
{\it strictly parabolic} if $\inf d_f=0$.

In the case of the hyperbolic spaces ${\bf H}^n$ this classification
is same as the standard one, and there are additional properties:
an axis of a hyperbolic isometry is unique, a parabolic
isometry has a fixed point in $X(\infty)$, which is unique.

\subsection{Parabolic isometries}
In general a parabolic isometry may not have any fixed
point in $X(\infty)$, but if $X$ is proper then
there is always at least one.
To see it we first quote a lemma from \cite{BH}.

\begin{lem}\label{displace}
Let $X$ be a proper $\CAT(0)$ space.
Let $\delta:X \to {\Bbb R}$ be a continuous,
convex function which does not attain its infimum.
Then there is at least one point $p \in X(\infty)$
such that if $g$ is an isometry of $X$ with $\delta(g(x))=\delta(x)$
for all $x \in X$, then $g$ fixes $p$ and leaves each horosphere
centered at $p$ invariant.
\end{lem}

For a proof we refer to 8.26 in Chap II.8 \cite{BH}.
This claim is stronger than the lemma in 8.26 there, however the
argument is same.
As a consequence we obtain the following proposition.

\begin{prop}\label{parabolicfix}
Let $X$ be a proper $\CAT(0)$ space. Suppose $f$ is a
parabolic isometry on $X$. Then there is a point
$p \in X(\infty)$ such that any isometry $g$ of $X$
with $gf=fg$, in particular, $f$, fixes $p$ and leaves each horosphere
centered
at $p$ invariant.
\end{prop}

Again this is similar but stronger than the proposition
in 8.25 in Chap II.8 \cite{BH}, where the claim
is only on $f$ itself and not about $g$ commuting with $f$.
The proof is identical too, but we put it since it is short.
We remark this claim is essentially same as
Lemma 7.3 (2) on p87 in \cite{BGS}, where
they consider only manifolds of non-positive curvature.

\begin{proof}
Consider the displacement function $d_f$ of $f$. It is a convex function
which does not attain its infimum.
Since $gf=fg$, $d_f$ is $g$-invariant.
Now we apply Lemma \ref{displace} putting $\delta=d_f$.
\end{proof}

A parabolic isometry of CAT(0) space may have
more than one fixed point in the ideal boundary.
For example, let $T$ be a tree and put $X=T \times {\bf H}^2$.
Let $g$ be a parabolic isometry on ${\bf H}^2$
and put $f=\text{id}_T \times g$, where id$_T$ is the
identity map on $T$.
Then $f$ is a parabolic isometry on $X$ and
the set of fixed points in $X(\infty)$, denoted by $X_f(\infty)$, is
the join of $T(\infty)$ and the (only)
fixed point of $g$ in ${\bf H}^2(\infty)$.
We remark that the conclusion of Proposition \ref{parabolicfix} may
not apply to all points $p$ in $X_f(\infty)$.

\section{Horospheres}
\subsection{Covering dimension}\label{covdim}
We recall the definition of the
covering dimension (cf.\ \cite{Ng:dim}).
Let ${\mathcal U}$ be a covering of a topological space.
The {\it order} of ${\mathcal U}$ at a point $p \in X$, ord$_p({\mathcal
U})$,
is the number (possibly $\infty$) of the members in ${\mathcal U}$ which
contain $p$.
The {\it order} of ${\mathcal U}$, ord(${\mathcal U}$),
is $\sup_{p \in X} \text{ord}_p({\mathcal U})$.

A covering ${\mathcal V}$ is a {\it refinement} of
a covering ${\mathcal U}$ if each member of ${\mathcal V}$
is contained in some member of ${\mathcal U}$.
If for any finite open covering ${\mathcal U}$ of a topological
space $X$ there exists an open covering ${\mathcal V}$ which is a
refinement of ${\mathcal U}$ such that ord(${\mathcal V}) \le n+1$, then
$X$ has covering dimension $\le n$.
$X$ has {\it covering dimension} $=n$ if the covering
dimension is $\le n$ but not $\le n-1$.

It is a theorem that in the case that $X$ is a metric space
one can replace ``any finite open covering" by ``any open covering"
in the definition of covering dimension (see II 5 in \cite{Ng:dim}).

\subsection{A $Z$-set in an ANR}
We state a theorem which is a consequence of
a general result by Bestvina-Mess, \cite{BM}.
They are interested mainly in the ideal
boundary of a word-hyperbolic
group.

Recall that a closed subset $Z$ in a compact
ANR $Y$ is a {\it $Z$-set}
if for every open set $U$ in $Y$ the inclusion
$U-Z \to U$ is a homotopy equivalence.
It is known (\cite{BM}) that
each of the following properties characterizes $Z$-sets:

\begin{enumerate}
\item[\rm(i)]
For every $\epsilon >0$ there is a map $Y \to Y-Z$,
which is $\epsilon$-close to the identity.
\item[\rm(ii)]
For every closed $A \subset Z$ there is a homotopy
$H:Y\times [0,1] \to Y$
such that $H_0=$ identity, $H_t|A=$ inclusion, and
$H_t(Y-A) \subset Y-Z$
for $t>0$.
\end{enumerate}
We first quote one result from their paper (Proposition 2.6 in \cite{BM}).

\begin{prop}\label{zset}
Suppose $Y$ is a finite-dimensional ${\rm ANR}$ and $Z \subset Y$
a $Z$-set. Then $\dim Z < \dim Y$, and hence
$\dim Z < \dim(Y-Z)$.
\end{prop}

We show the following.

\begin{thm}\label{bestvina}
Let $X$ be a proper $\CAT(0)$ space of dimension
$n < \infty$. Then
\begin{enumerate}
\item
A metric sphere $S(x,r) \subset X$
for $x \in X, r>0$
has dimension at most $n-1$.
\item
A horosphere centered at a point $p \in X(\infty)$
has dimension at most $n-1$.
\item
$\bar{X}$ is an AR of dimension $n$.
The ideal boundary $X(\infty)$ is a $Z$-set in $\bar{X}$,
and the dimension is at most $n-1$.
$\bar{X}$ is homotopy equivalent to $X$, hence
contractible.
\end{enumerate}
\end{thm}

\begin{proof}
First of all, for a metric space of finite dimension
to be an ANR is equivalent to being locally contractible.
Therefore a CAT(0) space $X$ of finite
dimension as well as any convex subset in $X$
is an ANR.

(1)\qua A (closed) metric ball $B=B(x,r)$ in $X$ is an ANR since
it is a convex subset.
And its boundary, namely, the metric sphere
$S=S(x,r) \subset B$ is a $Z$-set.
To see it one verifies (i) in the above.
Indeed for any $\epsilon >0$ one can construct,
using a unique geodesic from each point in $B$ to
the center of $B$, a map $B \to B-S$ which is $\epsilon$-close to
the identity.
It follows that $\dim S < \dim B \le \dim X =n$.

(2)\qua A closed horoball $B$ in $X$ is an ANR since it is convex,
and its boundary, namely, a horosphere, $H \subset B$
is a $Z$-set. As in (1) one can construct a map which satisfies
the property (i) using a unique geodesic from each point in $B$
to $p$.

(3)\qua This is more complicated than (1) and (2).
One can construct a homotopy to show
$X(\infty) \subset \bar{X}$
is a $Z$-set, but we will show it as one of the conclusions of
the proposition below.

We quote another proposition from \cite{BM}.

\begin{prop}\label{anr}
Suppose $Y$ is a compactum (i.e., compact, metrizable space) and $Z
\subset Y$
a closed subset such that
\begin{enumerate}
\item
${\rm int}(Z)=\emptyset$.
\item
$\dim Y =n < \infty$.
\item
For every $k=0,1,\cdots,n$, every point $z \in Z$,
and every neighborhood $U$ of $z$, there is a neighborhood
$V$ of $z$ such that every map $\alpha:S^k \to V-Z$
extends to $\tilde{\alpha}:B^{k+1} \to U-Z$.
\item
$Y-Z$ is an ${\rm ANR}$.
\end{enumerate}
Then $Y$ is an ${\rm ANR}$ and $Z \subset Y$ is a $Z$-set.
\end{prop}

Let $Y=\bar{X}$ and $Z=X(\infty)$.
$\bar{X}$ is a compact metrizable space because $X$ is proper.
We check (1)-(4).

(1) is clear.
To show (2), fix $x \in X$. Consider a system of metric balls
centered at $x$, $B(x,r), r>0$, and a system of metric
spheres, $S(x,r), r>0$. For $0<r \le s$, define
$p_{rs}:B(x,s) \to B(x,r)$ by
$p_{rs}(y)=y$ if $y \in B(x,r)$ and
$p_{rs}(y)=[x,y](r)$ if $y \in B(x,s)-B(x,r)$,
where $[x,y]$ is the geodesic from $x$ to $y$.
Remark in the latter case $p_{rs}(y) \in S(x,r)$,
so that it gives $p_{r,s}|S(x,s):S(x,s) \to S(x,r)$.
We obtained inverse systems $\{B(x,r), p_{rs}\}_{r,s>0}$
and $\{S(x,r), p_{rs}\}_{r,s>0}$.
Since $X$ is proper they are compact Hausdorff spaces.

The following is more or less direct from
the definition of the cone topology.
We leave a precise argument to readers.

\medskip
{\bf Claim}\qua
(a)\qua $\bar{X}$ with the cone topology is the inverse limit of the inverse
system $\{B(x,r), p_{rs}\}_{r,s>0}$.

\noindent
(b)\qua $X(\infty)$ with the cone topology is the inverse limit of
the inverse system $\{S(x,r), p_{rs}\}_{r,s>0}$.

\medskip
Since the dimension of the inverse limit of
compact Hausdorff spaces of
$\dim \le d$ is also at most $d$ (e.g., see
1.7 Cor, Ch8 in \cite{P}),
$\dim \bar{X} \le n$ and $\dim X(\infty) \le n-1$.
Therefore $\dim Y \le n$.

(3)\qua
There exist $x \in X, r>0, \epsilon >0$ such
that $U(x,r,\epsilon;z) \subset U$.
Let $V=U(x,r,\epsilon;z)$. Then $V-Z=V \cap X$ is
convex, therefore
one can just cone off $\alpha$ in $V-Z$
to extend it to $\tilde{\alpha}$, which is still in $V-Z$.

(4)\qua 
Since $Y-Z$ is $X$ itself, it is an ANR.

Now by Proposition \ref{anr},
$\bar{X}$ is a finite-dimensional (in fact $n$-dim)
ANR and $X(\infty)\subset \bar{X}$
is a $Z$-set.
Therefore $X \subset \bar{X}$ is homotopy equivalent,
so that $\bar{X}$ is contractible, hence AR.
\end{proof}


\subsection{The cohomological dimension of a parabolic
subgroup}\label{cohodim}

Let $\cd G$ denote the {\it cohomological dimension} of $G$.
It is defined (cf.\ p185 \cite{Br}) by
$$\cd G=\sup\{n:H^n(G,M)\not= 0 \text{ for some $G$-module $M$} \}.$$
It is known that $\cd G = \geomdim G$ if
$\cd G \ge 3$ (7.2 Corollary p205 \cite{Br}). Also
$G$ is free if and only if $\geomdim G=1$, as well as if and only if
$\cd G=1$.
It is known (6.7 Proposition on p202 \cite{Br}) that if $G$ is of type
FP, in particular,
if there is a finite $K(G,1)$, then
$$\cd G=\max\{n:H^n(G,{\Bbb Z}G) \not=0\}.$$
We use Theorem \ref{bestvina} to bound the cd of
a group leaving horospheres invariant.
We remark that horospheres are not connected, not simply connected
and not locally contractible in general.
But we can at least ``approximate" them with locally
finite simplicial complexes (maybe not connected),
equivariantly in terms of group actions.

\begin{prop}\label{graph}
Let $X$ be a proper $\CAT(0)$ space of
dimension $n < \infty$.
Let $p \in X(\infty)$ and $\{H_t \}_{t \in {\Bbb R}}$
be the family of horospheres centered at $p$.
Suppose a group $G$ is acting properly on $X$ by isometries
fixing $p \in X(\infty)$,
leaving each horosphere $H_t$ invariant.

If the action of $G$ on a horosphere $H_t$ is free
then there exists a locally finite simplicial complex of
dimension $\le n-1$, $L_t$,
such that $G$ acts freely on $L_t$ by simplicial isomorphisms
and there is a $G$-equivariant continuous map
$f_t:H_t \to L_t$.
\end{prop}

\begin{proof}
By Theorem \ref{bestvina}, the dimension of $H_t$
is at most $n-1$.
Since the
action of $G$ on $H_t$ is free and proper, we have a
covering $\pi : H_t \to H_t/G$. In this proof, let $B(x,r)$ denote the
open metric ball in $H_t$ centered at $x \in H_t$ and of radius $r > 0$
with respect to the metric on $X$. Since
the action of $G$ is proper,
for each $x \in H_t$ there exists
a number $\rho(x) > 0$ such that $B_x := B(x,\rho(x))$ does not
intersect $gB_x$ whenever $g \in G \setminus \{e\}$.
We may assume that for all $x\in X, g \in G$, $\rho(gx)=\rho(x)$,
so that $B_{gx}=gB_x$.
Since
$H_t/G$ has covering dimension $\le n-1$ and since $\{\pi(B_x)\}_{x
\in H_t}$ is an open covering of $H_t/G$, according to Corollary of
II.6 in \cite{Ng:dim} we can choose a locally finite refinement
$\{U\}$ of it with order $\le n$, i.e.,
for any compact set in $H_t/G$, there are
only finitely many $U$'s which intersect the compact
set,
and any distinct $n+1$ elements in
$\{U\}$ has empty intersection.

We take a lift $\hat U$ of each
$U$ in $H_t$. Then, by the properness of the action of $G$,
the family $\mathcal{U} := \{g\hat U\}_{U,g \in G}$ is a locally
finite open covering of $H_t$. Let us now prove that $\mathcal{U}$
has order $\le n$. In fact, if not, there would exist $n+1$ different open
sets $g_i\hat U_i$, $i=1,2,\cdots, n+1$, in $\mathcal{U}$ that contain a
common
point. We can find $x_i \in H_t$ and $g_i' \in G$ for each $i$ such
that $g_i\hat U_i$ is contained in $g_i' B_{x_i}$. Therefore, by
recalling the definition of $\rho$, $g_i\hat U_i$ does not intersect
$g \hat U_i$ for any $g \in G$ with $g \neq g_i$, and in particular, if
$i \neq j$, then $\hat U_j$ is different from $g\hat U_i$ for any $g
\in G$. Thus $\pi(\hat U_i)$ for $i = 1,2,\cdots,n+1$ are all different, but
have a common intersection point, which contradicts that the order of
$\{U\}$ is $\le n$. Thus $\mathcal{U}$ has order $\le n$.

Let $L_t$ be the nerve of $\mathcal{U}$, which is a
locally finite simplicial complex
of dimension $\le n-1$ since ord($\mathcal{U}) \le n$.
Denote its geometric
realization also by $L_t$. $L_t$ (and $H_t$) may not be connected.
Since $G$ acts freely on the set $\mathcal{U}$,
we have an induced free action of $G$ on $L_t$. We shall construct
a $G$-equivariant continuous map $f : H_t \to L_t$. Let $w_V(x) :=
d(x,H_t \setminus V)$ for $x \in H_t$ and $V \in \mathcal{U}$. Note
that $w_V$ is a continuous function with the property that $w_V > 0$
on $V$ and $w_V = 0$ on $H_t \setminus V$.
Using this set of functions as ``coordinates" we define
a map $f_t:H_t \to L_t$ as follows.
For a point $x \in H_t$,
let $V_1,V_2,\cdots,V_m \in \mathcal{U}$
be the open sets which contain $x$, so that $1 \le m \le n$.
By the definition of $L_t$ there is a simplex, $\sigma$, of dimension
$m-1$ which corresponds to the sets $V_i$'s.
Since $L_t$ is a simplicial complex we can put a system of linear
coordinates on all simplexes so that they are
canonical in terms of restricting them to faces
of each simplex.
We define $f_t(x) \in \sigma$ to be the point
with the coordinates
$$\left(\frac{w_{V_1}(x)}{w(x)},
\cdots,\frac{w_{V_m}(x)}{w(x)}\right),$$ where
$w(x)=\sum_{i=1}^m w_{V_i}(x)$.
By construction $f_t$ is a (well-defined) continuous and $G$-equivariant
map.
\end{proof}

\begin{prop}\label{cdG}
Let $X$ be a proper $\CAT(0)$ space of dimension $n < \infty$.
Suppose a group
$G$ is acting properly on $X$ by isometries.
Assume $G$ fixes a point $p \in X(\infty)$ and leaves each horosphere $H_t$
centered at $p$ invariant.
If $G$ has a finite $K(G,1)$-complex then
$\cd G \le n-1$.
In particular, if $n=2$ then $G$ is free.
\end{prop}

\begin{proof}
Since $G$ has a finite $K(G,1)$, $G$ is torsion-free, so that
the action on $X$ is free.
Fix a Busemann function $b_p$ associated to $p$, and
let $H_t$ and $B_t$ be the set of horospheres and horoballs defined by
$$H_t=\{ x \in X\,|\, b_p(x)=t\}, B_t=\{ x \in X\,|\, b_p(x) \le t\}.$$
Each of them is $G$-invariant.
Let $K$ be the universal covering of a finite $K(G,1)$ complex.
We claim that there exists $u$ such that
there is a $G$-equivariant continuous map $h:K \to H_u$.

We first construct a map $f$ from $K$ to $X$
such that the image of $K^0$ is in one horosphere.
Fix $s$ and give
a $G$-equivariant map $f:K^0 \to H_s$ arbitrarily.
Since $X$ is contractible, $K/G$ is
a finite complex and the action of $G$ on $K$ is free,
one can extend
this map continuously to a $G$-equivariant map $f:K \to X$.

We now ``project" the map $f$ so that the whole
image is contained in one horosphere (maybe different
from $H_s$).
By construction, there is a constant $C$ such that
$f(K) \subset B_{s+C} - B_{s-C}$.
Define a projection
$${\rm pr}:B_{s+C}-B_{s-C} \to H_{s-C}$$
by ${\rm pr}(x)=[x,p] \cap H_{s-C}$,
where $[x,p]$ is the unique geodesic from $x$ to $p$.
This projection is $G$-equivariant.
Let
$$h={\rm pr} \circ f :K \to H_{s-C}.$$
This is $G$-equivariant as well.
Put $u=s-C$. We obtained a desired map $h$.

By Proposition \ref{graph}, for each $t$
there is a $G$-equivariant continuous map
$f_t:H_t \to L_t$, where $L_t$ is a
locally finite simplicial complex of $\dim \le n-1$ on which $G$ acts freely
by simplicial isomorphisms.
Let's write $L_u$ as $H$ (because this is a substitution
to a horosphere) from now on and
consider the composition $k=f_u \circ h:K \to H$.
This is a continuous, $G$-equivariant map.

Since $K$ is contractible, and $H$ is a simplicial
complex on which $G$ acts freely by simplicial isomorphisms,
there exists a $G$-equivariant
homotopy inverse of $k$, $g: H \to K$, i.e.,
$g \circ k$ is $G$-equivariantly homotopic to the identity
of $K$.
We will show $H^i(G,{\Bbb Z}G)=0$ for $i > n-1$,
which implies $\cd G \le n-1$.

In general we have
$$H^*(G,M) = H^*(K(G,1);{\mathcal M}),$$
where $M$ is a $G$-module and ${\mathcal M}$ is
the local coefficient system on $K(G,1)$ associated
to $M$ (p59 \cite{Br}).
Therefore since $K/G$ is $K(G,1)$,
it suffices to show $H^i(K/G;{\Bbb Z}G)=0$ for $i > n-1$,
where we abuse a notation and ${\Bbb Z}G$ denotes
the local coefficient system on $K/G$ with
the natural action of $G=\pi_1(K/G)$ on ${\Bbb Z}G$.

Let $k,g$ denote the continuous maps induced by $k,g$
on $K/G,H/G$ as well.
Consider the induced homomorphisms
$g^*:H^i(K/G;{\Bbb Z}G) \to H^i(H/G;{\Bbb Z}G),
k^*:H^i(H/G;{\Bbb Z}G) \to H^i(K/G;{\Bbb Z}G).$
Then the composition $k^* \circ g^*$ is
an isomorphism of $H^i(K/G;{\Bbb Z}G)$
because the composition $g_*\circ k_*$
is the identity on $\pi_1(K/G)=G$, so that
the action of $G$ on ${\Bbb Z}G$ is the original one.
For $i > n-1$, clearly $H^i(H/G;{\Bbb Z}G)=0$
because $\dim (H/G)= \dim H \le n-1$, so that
$H^i(K/G;{\Bbb Z}G)=0$ as well.
\end{proof}

\section{Applications}

In this section we first apply Proposition \ref{cdG}
to free abelian groups and prove Proposition \ref{ss}.
Then we discuss the example by Brady-Crisp \cite{BC},
the example by Bridson \cite{Bri} and new examples
in connection to Problem \ref{problem}.
An interesting and new aspect about the new examples is that
although they do not act properly, nor even freely,
on any CAT(0) spaces by semi-simple
isometries at all, they act properly on some proper CAT(0) spaces
(indeed even CAT($-1$) spaces) if we allow
parabolic isometries. However, the minimal dimensions
of such CAT(0) spaces are strictly bigger
than the geometric dimensions. See Corollary \ref{heisenberg} for the
precise statement.

\subsection{Flat torus theorems}
The following is called the flat torus theorem.
See 7.1 in Chap II, \cite{BH} for a proof.

\begin{thm}[Flat torus theorem \cite{BH}]\label{flattorus}
Let $G$ be a free abelian group of rank $n$ acting properly on a
complete CAT(0) space
$X$ by semi-simple isometries.
Then there is a $G$-invariant, convex subspace in $X$ which is isometric
to the Euclidean space ${\bf E}^n$ of dimension $n$ such that the quotient
of the action by $G$ is an $n$-torus.
\end{thm}

\begin{remark}
In the above theorem, in fact, it suffices to assume
that the action is only K-proper. An argument is same
as the one in \cite{BH}. Of course, the freeness of the
action is not enough. Think of a free action of
${\Bbb Z}^2$ on ${\bf R}$ by isometries.
\end{remark}

The following is a main result of the paper.
We discuss several applications of this proposition in
the rest of the paper.

\begin{prop}\label{ss}
Let $X$ be a proper $\CAT(0)$ space.
Suppose a free abelian group of rank $n$, ${\Bbb Z}^n$,
is acting properly on $X$ by isometries.
Then $\dim X \ge n$.
If $\dim X =n$
then the action is by semi-simple isometries.

\end{prop}

\begin{proof}
First of all, the dimension of $X$ is at least $n$
since $X/{\Bbb Z}^n$ is a $K(\pi,1)$ space
and $\cd {\Bbb Z}^n =\geomdim {\Bbb Z}^n =n$.

Suppose there is a parabolic isometry $f \in {\Bbb Z}^n$ on
$X$.
Then since $X$ is proper, by Proposition \ref{parabolicfix} there exists
$p \in X(\infty)$
such that ${\Bbb Z}^n$ fixes $p$ and leaves each horosphere
centered at $p$ invariant.
Since ${\Bbb Z}^n$ has a finite $K(\pi,1)$ it follows
from Proposition \ref{cdG}
that $n=\cd {\Bbb Z}^n \le \dim X-1$.
Therefore $n+1 \le \dim X$.
\end{proof}

We obtain another torus theorem.

\begin{thm}[Flat torus theorem]
Suppose $X$ is a proper geodesic space
of dimension $n\ge 2$ whose fundamental
group is ${\Bbb Z}^n$.
If the universal covering of $X$ is $\CAT(0)$ then
$X$ contains a totally geodesic flat $n$-torus
which is a deformation retract.
\end{thm}

\begin{proof}
Consider the proper
isometric action of $G={\Bbb Z}^n$ on the universal cover $Y$ of $X$.
By Proposition \ref{ss} the action is by semi-simple isometries.
Then by Theorem \ref{flattorus}, there is $P={\bf E}^n \subset Y$
which is invariant by $G$.
The flat $n$-torus $P/G \subset X$ is a deformation retract.
The retraction is given
by the projection of $Y$ to $P$, which is $G$-equivariant.
\end{proof}

\subsection{Nilpotent groups, Heisenberg groups and complex hyperbolic
manifolds}

\begin{thm}\label{nilpotent}
Let $G$ be a torsion-free,
nilpotent group.
\begin{enumerate}
\item
If $G$ acts freely
on a $\CAT(0)$ space by semi-simple isometries,
then $G$ is abelian.
\item
Suppose $G$ is finitely generated, not abelian,
and acting properly
on a proper $\CAT(0)$ space $X$ by isometries.
Then
$\cd G \le \dim X-1$.
\item
Suppose $G$ is finitely generated, and
acting properly on a proper $\CAT(0)$ space $X$.
If $\cd G=\dim X$ then $G$ is abelian and
the action is by semi-simple isometries.
\end{enumerate}
\end{thm}

\begin{proof}
(1)\qua This claim is stated for manifolds of non-positive curvature
in 7.4 Lemma (3) in \cite{BGS}. Our argument is essentially same.
Suppose $G$ is not abelian. Let $C(G)$ be the center.
There exists a non-trivial element $g \in C(G) \cap [G,G]$.
Then $g$ generates ${\Bbb Z}$ in $G$.
If $G$ is acting on a CAT(0) space $X$ freely
by isometries then $g$ has to be parabolic.
This is because if $g$ was
hyperbolic then the union of the axes of $g$ in $X$
would be isometric to $Y \times {\bf R}$, which is $G$-invariant.
Moreover each $h \in G$ acts on $Y \times {\bf R}$
as a product of isometries on $Y$ and ${\bf R}$.
The induced isometric actions of $G$ on ${\bf R}$
is by translations because $g$ is a non-trivial
translation and in the center.
But since $g \in [G,G]$ the action of $g$ on ${\bf R}$
must be trivial, which is a contradiction.

(2)\qua Note that the action is free because
$G$ is torsion-free. Take $g$ as in (1), so that
$g$ is parabolic.
Since $g \in C(G)$, by Proposition \ref{parabolicfix}, there is a point
$ p \in X(\infty)$
such that $G(p)=p$ and each horosphere $H$ at $p$ is $G$-invariant.
It is a theorem of Malcev that there exists a simply connected
nilpotent Lie group $L$ which contains $G$ as a co-compact
lattice, so that $L/G$ is a closed manifold of finite dimension
with $\pi_1 \simeq G$. It follows that $G$ has a finite $K(\pi,1)$-complex
(for example take the nerve of a covering of $L/G$).
Therefore by Proposition \ref{cdG}, $\cd G \le \dim H \le \dim X-1$.

(3)\qua By (2), $G$ is abelian, and its rank is $\cd G$.
Then we already know that the action
has to be by semi-simple isometries by Proposition \ref{ss}.
\end{proof}


Let ${\bf H}_{\Bbb C}^n$ denote
the complex hyperbolic space of complex-dimension
$n$ (cf.\ \cite{CG}).
${\bf H}_{\Bbb C}^n$ is a complete, $2n$-dimensional
Riemannian manifold of sectional curvature
bounded between $-1$ and $-1/4$ .
It is a proper CAT($-1$) space of dimension $2n$.

\begin{cor}\label{heisenberg}
For an integer $n \ge 1$, let
$G_n$ be the finitely presented group defined by a set of generators
$\{ a_1,b_1, \cdots, a_n,b_n,c \} $ and a set of
relators as follows: all commutators in generators
$=1$ except for $[a_i,b_i]=c$ for all $i$.

Then $G_n$ is a finitely generated, torsion-free,
nilpotent group with the following properties.

\begin{enumerate}
\item
$\geomdim G_n=2n+1$.
\item
$G_n$ does not act properly on any proper $\CAT(0)$ space
of dimension $2n+1$ by isometries.
\item
$G_n$ acts properly on ${\bf H}_{\Bbb C}^{n+1}$, which is a proper
$\CAT(-1)$ space of dimension $2n+2$,
by isometries such that all
non-trivial elements are strictly parabolic.
\item
If $G_n$ acts freely on a $\CAT(0)$ space
by isometries,
then $c$ has to be parabolic. In particular,
$G_n$ does not act freely, hence not properly,
on a $\CAT(0)$ space by semi-simple isometries,
so that $G_n$ is not a $\CAT(0)$ group.
\end{enumerate}
\end{cor}

\begin{proof}
It is known that
$G_n$ is a finitely generated, torsion free, nilpotent
group such that $[G_n,G_n]=\langle c \rangle =C(G_n)\simeq{\Bbb Z}$,
and $G_n/[G_n,G_n] \simeq {\Bbb Z}^{2n}$.
They are called {\it (discrete) Heisenberg groups}.

(1)\qua See p186 in \cite{Br}.
It is known that the {\it rank} (or {\it Hirsch number})
of $G_n$ is equal to $\cd G_n$.
And the rank of $G_n$ is $2n+1$
because $G_n$ has a central series $G_n > [G_n,G_n]\simeq {\Bbb Z}
>1$ with $G_n/[G_n,G_n]\simeq {\Bbb Z}^{2n}$.

Remark that if $n=1$ then $G$ is the fundamental group of
a closed $3$-manifold, which is a non-trivial circle bundle
over a torus.

(2)\qua Combine (1) and Theorem \ref{nilpotent} (2).

(3)\qua
It is known that one can embed
$G$ in $\text{Isom}({\bf H}_{\Bbb C}^{n+1})$
as a discrete subgroup, so that
$G$ acts properly on ${\bf H}_{\Bbb C}^{n+1}$
by strictly parabolic isometries.
To see it, fix a point in the ideal boundary and
let $P_n < \text{Isom}({\bf H}_{\Bbb C}^{n+1})$ be the
subgroup of all strictly parabolic isometries (and the trivial
element) fixing
this point. Then $P_n$ is a nilpotent Lie group which is
a central extension of ${\Bbb R}$ by ${\Bbb C}^n$.
$P_n$ is called the Heisenberg Lie group.
If we take ``integer points" in $P_n$ we obtain
$G_n < P_n$ as a discrete subgroup.
One finds details on this paragraph in Proposition 4.1.1 in \cite{CG}.

(4)\qua Suppose $c$ is semi-simple on $X$. Let $Y \subset X$
be the union of axes of $c$. $Y$ is isometric
to $A \times {\bf R}$.
But since $c$ is in the center,
as we saw in the proof of Theorem \ref{nilpotent} (1),
$c$ has to act trivially on ${\bf R}$ because $c=[a_1,b_1]$,
a contradiction.
\end{proof}

\begin{question}
Does $G_n$ act properly on a CAT(0) space of dimension $2n+1$
by isometries? (No if the space is proper).
\end{question}

\begin{cor}\label{complex}
Let $G$ be the fundamental group of a non-compact, complete,
complex hyperbolic manifold $M$
of complex-dimension $n \ge 2$. If the volume of $M$ is finite, then

\begin{enumerate}
\item
$\cd G=\geomdim G=2n-1$.
\item
$G$ does not act properly on
any proper $(2n-1)$-dimensional $\CAT(0)$ space
by isometries.
\item
$G$ acts properly on ${\bf H}_{\Bbb C}^{n}$, which is a proper
$2n$-dimensional $\CAT(-1)$ space, by isometries.
\item
$G$ does not act freely, hence not properly, on any $\CAT(0)$ space
by semi-simple isometries,
so that $G$ is not a $\CAT(0)$ group.

\end{enumerate}
\end{cor}

\begin{proof}
(1)\qua First of all $\dim M=2n$.
Since $M$ is not compact, by
8.1 Proposition \cite{Br}, $\cd G < 2n$.
Because the volume of $M$ is finite, it has
cusps. Let $N$ be a cross-section
of one of the cusps, and $P$ its fundamental
group. $P$ is torsion-free. $N$ is a closed, infra-nil manifold
(i.e., finitely covered by a nil manifold)
of dimension $2n-1$.
The universal cover of $N$ is diffeomorphic
to ${\bf R}^{2n-1}$ since
it is a horosphere of ${\bf H}_{\Bbb C}^{n}$.
Thus $\cd P=2n-1$, so that $\cd G=2n-1$
as well because $P<G$ and $\cd G < 2n$.
In general $\cd G=\geomdim G$ if $\cd G \ge 3$.

(2)\qua
Because $N$ is an infra-nil manifold,
$P$ is virtually nilpotent. But $P$ is not virtually abelian.
Let $P'<P$ be a nilpotent group of finite index. Then $\cd P' =2n-1$ and
$P'$ is finitely generated. To argue by contradiction suppose $G$ acts on a
proper $\CAT (0)$ space $X$ of dimension $2n-1$ properly by isometries.
We apply Theorem \ref{nilpotent} (3) to the action of $P'<G$ on $X$ and
conclude that $P'$ is abelian. Impossible since $P$ is not virtually
abelian.


(3)\qua The universal cover of $M$ is isometric to ${\bf H}_{\Bbb C}^n$,
which is
a proper CAT($-1$) space of dimension $2n$.
The action of $G$ on the universal cover of $M$ is a desired one.

(4)\qua This is similar to (2).
To argue by contradiction suppose we had such action. The action is free.
Let $P'<G$ be as in (2).
Apply Theorem \ref{nilpotent} (1) to the action of $P'$ and get a
contradiction.
\end{proof}

A parallel result holds for quaternionic hyperbolic spaces as well.
Let ${\bf H}_{\Bbb H}^n$ denote
the quoternionic hyperbolic space of quoternionic-dimension
$n$ (cf.\ \cite{CG} for a unified treatment of
real, complex and quoternionic hyperbolic spaces).
${\bf H}_{\Bbb H}^n$ is a complete, $4n$-dimensional
Riemannian manifold of sectional curvature
bounded between $-1$ and $-1/4$ .
It is a proper CAT($-1$) space of dimension $4n$.

\begin{cor}\label{quaternion}
Let $G$ be the fundamental group of a non-compact, complete,
quoternionic hyperbolic manifold $M$
of quoternionic dimension $n \ge 2$. If the volume of $M$ is finite, then

\begin{enumerate}
\item
$\cd G=\geomdim G=4n-1$.
\item
$G$ does not act properly on
any proper $(4n-1)$-dimensional $\CAT(0)$ space
by isometries.
\item
$G$ acts properly on ${\bf H}_{\Bbb H}^{n}$, which is a proper
$4n$-dimensional $\CAT(-1)$ space, by isometries.
\item
$G$ does not act freely, hence not properly, on any $\CAT(0)$ space
by semi-simple isometries,
so that $G$ is not a $\CAT(0)$ group.

\end{enumerate}
\end{cor}

\begin{proof}
An argument is nearly identical to the proof of
Cor \ref{complex}. Just replace ${\bf H}_{\Bbb C}^{n}$
by ${\bf H}_{\Bbb H}^{n}$, and $2n$ by $4n$ appropriately.
\end{proof}

Interestingly, there are different
phenomena in the real hyperbolic case.
Let $M$ be the mapping torus
of a homeomorphism, $f$, of the
once puncture torus.
$G=\pi_1(M)$ is an extension of $F_2$ by ${\Bbb Z}$
such that ${\Bbb Z}$ acts on $F_2$ by an automorphism,
$\phi$.
$$ 1 \to F_2 \to G \to {\Bbb Z} \to 0.$$
T.Brady, \cite{TBr}, showed the following.

\begin{thm}\label{brady}
Suppose $G$ does not contain the
direct product $F_2 \times {\Bbb Z}$ as a
subgroup of finite index.
Then $G$ acts properly and co-compactly on some proper CAT(0) space of
dimension $2$ by isometries.
\end{thm}

He takes a $2$-dimensional spine, $X$, of $M$, which is
a $2$-complex, and put a metric which
makes $X$ into a piecewise Euclidean $2$-complex.
He shows that the universal cover of $X$ is a CAT(0) space.
Since $X$ is a retraction of $M$, the fundamental
group is $G$, so that its action on the universal
cover is a desired one.
Note that the action is by semi-simple isometries.

On the other hand, it is known that if (and only if)
$f$ is a pseudo Anosov map
then $M$ admits
a complete hyperbolic metric such that the volume
is finite, with one cusp.
A general result is due to Thurston (cf.\ Thm15.18, \cite{Kapo}).
There is a construction of a hyperbolic metric
using ideal hyperbolic tetrahedra, which is originated by J{\o}rgensen (cf.\
\cite{TBr}).
If $M$ is a hyperbolic manifold, then
its universal cover of $M$ is isometric to ${\bf H}^3$,
on which $G$ acts properly by isometries. Note that
a cusp subgroup is isomorphic to ${\Bbb Z}^2$,
and acts on ${\bf H}^3$ by parabolic isometries.


\subsection{Artin groups}

The Artin group $A(m,n,p)$ is defined by the following
presentation.
$$A(m,n,p)=\langle a,b,c\,|\,(a,b)_m=(b,a)_m, (b,c)_n=(c,b)_n,
(a,c)_p=(c,a)_p \rangle ,$$
where $(a,b)_m$ is the alternating product of length $m$ of $a$'s and $b$'s
starting with $a$.
Brady and Crisp (\cite{BC} Thm3.1) showed the following theorem.

\begin{thm}[Brady-Crisp]\label{BC}
Let $A=A(m,n,2)$.
Then we have the following.
\begin{enumerate}
\item
Let $m,n \ge 3$ be odd integers. Then
$A$ does not act properly on any complete $\CAT(0)$ space of
dimension $2$ by semi-simple isometries.
\item
With finitely many exceptions, for each $A$,
there is a proper $3$-dimensional
$\CAT(0)$ space $X$ on which $A$ acts
properly
by isometries such that $X/A$
is compact, hence the action is by semi-simple isometries,
and in particular, $A$ is a $\CAT(0)$ group.

\end{enumerate}
\end{thm}

\begin{remark}
Theorem \ref{BC} (1) remains true if we weaken
``properly" to ``K-properly" in the statement.
In their argument \cite{BC}, a certain properness of the action
is used in Th 1.3 (flat torus theorem), Lemma 1.4,
Prop 2.2, where the K-properness is enough. Mostly
it is about the properness of an action on a line.
Combining those, they obtain flats for
several ${\Bbb Z}^2$ subgroups in the beginning of
\S 3, where the theorem is proven, then they do not
use any properness of the action any more.

\end{remark}

It is known that $A(m,n,p)$ has a $2$-dimensional finite $K(A,1)$-complex
if $\frac{1}{m}+\frac{1}{n}+\frac{1}{p} \le 1$ (\cite{CD}),
so that $\geomdim A=2$ if $\frac{1}{m}+\frac{1}{n} \le \frac{1}{2}$.
Therefore there are infinitely many $A=A(m,n,2)$ such
that $\geomdim A=2$ but $\dimss A=3$.

Brady-Crisp asked if the conclusion of Theorem \ref{BC} (1) remains
valid if one drops
the assumption that the action is by semi-simple isometries.
Using Proposition \ref{ss} we give an answer to this question
in the case that $X$ is proper and of dimension $2$.

\begin{thm}\label{BC2}
Let $m,n \ge 3$ be odd integers.
The Artin group $A=A(m,n,2)$ does not act properly
on any proper $\CAT(0)$ space
of dimension $2$
by isometries.
\end{thm}

\begin{proof}
Brady-Crisp showed that if $A$ is acting on a
complete $2$-dimensional CAT(0) space $X$ properly
by isometries then $A$ contains a parabolic isometry, $g$.
In fact their argument shows an additional property on $g$ that
we may assume that the element $g$ is in a subgroup isomorphic to ${\Bbb
Z}^2$.

To prove our theorem, let's assume that such action of $A$
exists. Then as we have just said
there is a subgroup $G < A$ isomorphic
to ${\Bbb Z}^2$ which contains a non semi-simple isometry.
But this is impossible by Proposition \ref{ss}.
\end{proof}

Theorem \ref{BC} (2) and Theorem \ref{BC2} imply
Theorem \ref{cat0}.

\begin{question}
Does $A$ in Theorem \ref{BC} act properly on a complete CAT(0) space of
dimension $2$
by isometries?
(No if the space is proper).
\end{question}

\subsection{An example by Bridson}
Let $B$ be a group with the following presentation.
$$B=\langle a,b,\gamma,s,t\,|\,\gamma a \gamma^{-1}=a^{-1},
\gamma b \gamma^{-1}=b^{-1}, sas^{-1}=[a,b]=tbt^{-1} \rangle.$$
A geodesic space is {\it locally CAT$(0)$}
if its universal cover is CAT(0).
Bridson \cite{Bri} showed the following theorem.

\begin{thm}[Bridson]\label{bridson}$\phantom{99}$
\begin{enumerate}
\item
$\geomdim B=2$.
\item
$B$ does not act properly
on any complete $2$-dimensional $\CAT(0)$ space
by semi-simple isometries.
\item
$B$ is the fundamental group of a compact, locally $\CAT(0)$,
$3$-dimensional
cubed complex. In particular, $B$ is a $\CAT(0)$ group.
\end{enumerate}
\end{thm}

\begin{remark}
As in Theorem \ref{BC}, one may wonder if one can weaken
the condition ``properly" to ``K-properly" in the
point (2) of this theorem. We do not know the answer.
The only place he may really need the properness, not
just the K-properness of an action in his argument is
Prop 1.2, where he finds
a minimal tree of an action by $F_2$ on some ${\bf R}$-tree.
He uses a fact that if the action if proper, then the minimal
tree is simplicial.
\end{remark}

In (3), $B$ acts on the universal cover of the cubed
complex by semi-simple isometries, so that taking (2)
into account we conclude that $\dimss B=3$.
We show the following.
\begin{thm}\label{bridson2}
$B$ does not act properly
on any proper CAT(0) space of dimension $2$
by isometries.
\end{thm}

\begin{proof}
Suppose it did.
In the proof of Theorem \ref{bridson} (2), (cf.\ Proposition 3.1 in
\cite{Bri}), he
in fact showed that if $B$ acts properly by isometries
on a complete CAT(0) space of dimension $2$ then $\gamma^2$
must be parabolic.
On the other hand, the subgroup generated by $\gamma^2$ and $a$ is
isomorphic
to ${\Bbb Z}^2$, so that each element in this subgroup
is semi-simple by Proposition \ref{ss}, a contradiction.
\end{proof}

We remark that Theorem \ref{bridson} (3) and Theorem \ref{bridson2} imply
Theorem \ref{cat0} as well.

\begin{question}
Does $B$ act on a complete CAT(0) space of
dimension $2$
properly by isometries?
(No if the space is proper).
\end{question}

\subsection{Solvable Baumslag-Solitar groups}
A Baumslag-Solitar group is defined for non zero
integers $n,m$ by
$$BS(n,m)=\langle a,b\,|\,a b^n a^{-1}=b^m \rangle .$$
Suppose $|m| \ge 2$ in the following discussion.

(1)\qua $BS(1,m)$ is torsion-free and solvable, and its geometric
dimension is $2$.
One can construct a finite CW-complex of dimension $2$
which is a $K(\pi,1)$-space for $BS(1,m)$ as follows.
Take an annulus and glue one of its two
boundary components to the other one such
that it wraps $m$-times.

(2)\qua Suppose $BS(1,m)$ is acting
on a complete CAT(0) space freely
by isometries.
Then $b$ is parabolic, because if not,
then the relation
$a b a^{-1}=b^m$ would force the translation length of
$b$, $l=\min d_b$, to satisfy $l=|m|l$, so that $l=0$.
This is impossible since $b$ has an infinite order.

(3)\qua $BS(1,m)$ acts on ${\bf H}^2$ by isometries freely.
Define a presentation $\phi:BS(1,m) \to SL(2,{\Bbb R})$
by
$$
\phi(a)=\left( \begin{array}{cc}
\sqrt{|m|} & 0 \\
0 & \sqrt{|m|}^{-1}
\end{array}
\right),
\phi(b)=
\left( \begin{array}{cc}
1 & 1 \\
0 & 1
\end{array}
\right).
$$
Since $PSL(2,{\Bbb R})=\text{Isom}({\bf H}^2)$,
$\phi$ gives a desired action, such that $a$ acts as a hyperbolic
isometry and $b$ is parabolic.

But $G$ does not act on ${\bf H}^n$ properly
by isometries. To see this, suppose it did. Then $b$ has to be parabolic
by (2).
Let $\alpha \in {\bf H}^n(\infty)$ be a point fixed by $b$.
Since $\alpha$ is a unique fixed point of $b$ and $b^m$ as well,
$a$ has to fix $\alpha$, so that $G$ fixes $\alpha$.
$G$ acts on the set of horospheres centered at $\alpha$.
Using a Busemann function associated to $\alpha$ we obtain
an action of $G$ on ${\bf R}$ by translations.
$b$ is in the kernel of the action.
The action of $G$ factors through the abelianization of $G$.
Thus the kernel of the abelianization, which is isomorphic to ${\Bbb Z}[
\frac{1}{m}]$,
leaves each of the horospheres invariant, and the action
on each horosphere is proper and by isometries.
Each horosphere is isometric to ${\bf E}^{n-1}$, so that
only a virtually free abelian group of finite rank
can act on it properly by isometries
by the Bieberbach theorem. But
${\Bbb Z}[ \frac{1}{m}]$ is not a virtually free abelian group of finite
rank,
a contradiction.

(4)\qua
$BS(1,m)$ is an HNN-extension of ${\Bbb Z}$ as follows:
$BS(1,m)={\Bbb Z}*_{{\Bbb Z}, f}$ such that the injective
homomorphism $f$ is given by $z \mapsto mz, z \in {\Bbb Z}$.
The other injective homomorphism is the identity map.
On the Bass-Serre tree, $T$, the element $b$ is elliptic
and $a$ is hyperbolic. We made $T$ into a CAT(0) space
such that each edge has length one. The tree $T$
is a regular tree of index $m+1$.
Note that for any vertex $v \in T$,
the stabilizer is generated by some conjugate of $b$.

Set $X={\bf H}^2 \times T$. This is
a proper CAT(0) space of
dimension three.
Using the action on ${\bf H}^2$
given in (3) and the Bass-Serre action on $T$,
we obtain the product isometric action of $BS(1,m)$ on $X$.
The action is proper.
To see this, let $x=(y,v) \in X$ be any point.
We may assume that $v \in T$ is a vertex.
Moreover, since the action of $BS(1,m)$ on $T$
is transitive on the set of vertices, we may assume
that the stabilizer subgroup of $v$ is generated
by $b$.
Set $r=1/2$. Then there are only finitely many
elements $g \in G$ with $d(x,gx) \le r$.
Indeed, if the inequality holds for $g$, then
$d(v,gv) \le 1/2$ on $T$, so that $g$
has to fix $v$. So, $g \in \langle b \rangle $.
On the other hand, $d(y,gy) \le 1/2$ on ${\bf H}^2$, but
there are only finitely many elements of the form
$b^n$ with $d(y,b^ny) \le 1/2$.


(5)\qua
Suppose $BS(1,m)$ acts on a proper
CAT(0) space $X$ properly by
isometries.
We recall a theorem by Adams-Ballmann.
Remark that ${\bf E}^0$ means a point.

\begin{thm}[Adams-Ballmann \cite{AB}]\label{AB}
Let $X$ be a proper $\CAT(0)$ space. Suppose
an amenable group $G$ acts on $X$ by isometries, then
either
\begin{enumerate}
\item[\rm(i)]
there is a point $\alpha \in X(\infty)$ fixed by $G$, or
\item[\rm(ii)]
there is a $G$-invariant convex subspace in $X$ which
is isometric to a Euclidean space ${\bf E}^n,(n \ge 0)$.
\end{enumerate}
\end{thm}

Since $G=BS(1,m)$ is solvable, hence amenable, by this theorem
there is a point $\alpha \in X(\infty)$ fixed by $G$.
This is because the case (ii) does not occur
by the Bieberbach theorem.
$G$ acts on the set of horospheres centered
at $\alpha$. As in (3), this gives a homomorphism $h$ from $G$ to ${\Bbb
R}$.
Clearly $h(b)=0$.

Now suppose $\dim X=2$ (we do not know if such actions exist. See the
questions below).
Then the homomorphism $h$ has to be non-trivial, since otherwise
$G$ leaves a horosphere invariant, which would imply
that $G$ is free by Proposition \ref{cdG}.
It follows that $h(a)\not=0$.
Therefore the kernel of $h$ is
isomorphic to ${\Bbb Z}[\frac{1}{m}]$.
This subgroup is acting on a horosphere properly
by isometries.

We record our discussion on $BS(1,m)=\langle a,b \,|\, aba^{-1}=b^m \rangle$
as follows.

\begin{prop}\label{bs}
Let $G=BS(1,m)$ such that $|m|\ge 2$.
\begin{enumerate}
\item
$\geomdim G=2$.
\item
If $G$ acts freely
on a $\CAT(0)$ space by isometries then
$b$ is parabolic. In particular, $G$ does
not act freely, hence not properly, on any $\CAT(0)$ space by
semi-simple isometries,
so that $G$ is not a $\CAT(0)$ group.
\item
$G$ acts freely on ${\bf H}^2$ by isometries.
But $G$ does not act properly on ${\bf H}^n$ for any $n \ge 1$
by isometries.
\item
$G$ acts properly on ${\bf H}^2 \times T$
by isometries where $T$ is a regular
tree of index $m+1$.
${\bf H}^2 \times T$
is a proper CAT(0) space of dimension $3$.


\item
If $G$ acts properly on a proper $\CAT(0)$ space $X$
by isometries, then
there exists a point $\alpha \in X(\infty)$ which is fixed by $G$.
The action of $G$ on the set of horospheres centered at $\alpha$
gives a homomorphism from $G$ to ${\Bbb R}$.
$b$ is in the kernel.
If $\dim X=2$ then the kernel of the homomorphism
is isomorphic to ${\Bbb Z}[\frac{1}{m}]$.

\end{enumerate}
\end{prop}

\begin{question}\label{bs.question}
Does $BS(1,m)$ act properly
on a CAT(0) space $X$ of dimension $2$?
\end{question}
If the answer is yes, then ${\Bbb Z}[\frac{1}{m}]$ acts on a horosphere
of $X$.
We do not know the answer even for proper spaces.
See Prop \ref{prop.bs.vis} for the case of visible CAT(0) spaces.
In connection to Theorem \ref{AB} we ask another question.

\begin{question}
Suppose an amenable group $G$ acts properly on a proper,
finite-dimensional
CAT(0) space $X$ by isometries.
Suppose $\cd G=\dim X$. Then is $G$ virtually abelian?
\end{question}

We know the answer is yes if $G$ is a finitely generated, nilpotent group
(Theorem \ref{nilpotent}). For example,
we do not know the answer for $BS(1,2)$, which is solvable, hence amenable.

\subsection{Torus bundle over a circle}
We discuss another example, which is
suggested by A.Casson.
Let $S$ be the group given by the following presentation.
$$S=\langle a,b,c\,|\,ab=ba,cac^{-1}=a^2b,cbc^{-1}=ab \rangle .$$
$S$ is solvable (but not virtually nilpotent), torsion free, and of
cohomological
dimension $3$ because it is the fundamental group
of a closed three-manifold, $M$, which is a torus bundle over a circle.
$M$ is a $K(\pi,1)$ space of $S$.
The subgroup generated by $a,b$, which is isomorphic to $\Z^2$,
is the fundamental group of the torus fiber, and
the base circle gives the element $c$.

$S$ is an HNN extension as $\Z ^2 *_{\Z ^2,f}$ such that
the self monomorphism $f$ of $\Z^2$ is the linear
map given by a matrix
$A=\left( \begin{array}{cc}
2 & 1 \\
1 & 1
\end{array}
\right) \in SL(2,\Z)
$.
The other monomorphism to define
the HNN extension is the identity.

$S$ is a semi-direct product
with the action of $\Z$ on $\Z^2$ given by $A$:
$$0 \to \Z ^2 \to S \to \Z \to 0.$$
Using the same matrix $A$, we produce a semi-direct product of
$\R^2$ and $\Z$ as follows. $S$ is a subgroup of $G$.
$$0 \to \R ^2 \to G \to \Z \to 0.$$
The group multiplication in $G$ for $(x,n), (y,m) \in {\Bbb R}^2 \times
\Z$ is
$$(x,n)(y,m)=(A^mx+y,n+m).$$
To define a monomorphism
$\sigma: G \to SL(2,{\Bbb R}) \times
SL(2,{\Bbb R})$, let
$p,1/p$ be the (positive, real) eigenvalues
of $A$, with eigenvectors $\alpha,\beta \in \R^2$.
$p=(3+\sqrt{5})/2$.
For $x \in \R^2$, write it uniquely as
$x=a_x\alpha + b_x \beta, a_x,b_x \in {\Bbb R}$, and define
$$\sigma(x,0)=
\left(
\left( \begin{array}{cc}
1 & a_x \\
0 & 1
\end{array}
\right),
\left( \begin{array}{cc}
1 & b_x \\
0 & 1
\end{array}
\right)
\right).$$
For $n \in \Z$, define
$$
\sigma(0,n) =
\left(
\left( \begin{array}{cc}
1/\sqrt p & 0 \\
0 & \sqrt p
\end{array}
\right)^n,
\left( \begin{array}{cc}
\sqrt p & 0 \\
0 & 1/\sqrt p
\end{array}
\right)^n
\right)
.$$
Since $(x,n) \in {\Bbb R}^2 \times {\Bbb Z}$
is written as a product
$(x,n)=(0,n)(x,0)$ in $G$, we define
$$\sigma(x,n)=\sigma(0,n)\sigma(x,0)
=
\left(
\left( \begin{array}{cc}
1/\sqrt p^n & a_x/\sqrt p^n \\
0 & \sqrt p ^n
\end{array}
\right),
\left( \begin{array}{cc}
\sqrt p^n & b_x \sqrt p^n \\

0 & 1/\sqrt p ^n
\end{array}
\right)
\right).
$$
Clearly $\sigma$ is a one to one map on $G$.
The map $\sigma$ is a homomorphism, namely
$\sigma(x,n) \sigma(y,m)=\sigma(A^mx +y,n+m)$.
One can check this using
$$A^m x= a_x p^m \alpha + \frac{b_x}{p^m} \beta.$$
We have used that $\alpha,\beta$ are eigenvectors of
$A$.

If we restrict $\sigma$ to $S$, the image
$\sigma(S)$ is a discrete subgroup in
$SL(2,{\Bbb R})
\times
SL(2,{\Bbb R})$, because
if a sequence $\sigma(x_i,n_i)$ converges,
then $n_i$ eventually becomes constant,
so that $a_{x_i},b_{x_i}$ become
constant as well since they are discrete
in ${\Bbb R}^2$.
$PSL(2,{\Bbb R})$ is by definition
$SL(2,{\Bbb R})/{\pm 1}$, which is
isomorphic to $\Isom_+({\bf H}^2)$,
the group of orientation preserving
isometries of ${\bf H}^2$.
Naturally we obtain a homomorphism,
which we also denote by $\sigma$,
$\sigma:S \to PSL(2,{\Bbb R})
\times
PSL(2,{\Bbb R})$.
This is also a monomorphism, so that
we obtain a proper,
isometric action of $S$ on $X={\bf H}^2 \times {\bf H}^2$.
The action is free because $S$ is torsion-free.
Therefore we obtain a non-compact, complete
$4$-dimensional Riemannian manifold
$N=({\bf H}^2 \times {\bf H}^2)/S$.

The manifold $N$ is homeomorphic to $M \times {\bf R}$.
To see this, let $\gamma$ be the geodesic in ${\bf H}^2$
which is the positive part of the $y$-axis in the
upper half plane model of ${\bf H}^2$.
Suppose $\gamma(t)$ is parametrized by arc length, upwards, i.e.,
the $y$-coordinate increases as $t$ increases.
Note that $\gamma$ is the common unique axis of hyperbolic
isometries $\sigma(0,n), n\not =0$, and
$v=\gamma(\infty)$ is the common unique fixed point
of $\sigma(x,0), x \not = 0$, so that
$v$ is also the unique fixed point of the action of $S$
on ${\bf H}^2(\infty)$.
The ideal boundary of ${\bf H}^2 \times {\bf H}^2$
is the spherical join of two copies of ${\bf H}^2(\infty)$.
The spherical join of two $v$'s in $({\bf H}^2 \times {\bf H}^2)(\infty)$
is the set of fixed point of $S$, which is a segment of length $\pi/2$.
Let $w$ be the midpoint of this segment, which is
$(\gamma \times \gamma)(\infty)$. Then each horosphere in $X={\bf H}^2
\times {\bf H}^2$,
$H_s, s \in {\Bbb R}$, centered at $w$ is invariant by $S$.
One can check this by the definition of $\sigma$.
We remark that if we take a point different from $w$ in the segment of
the fixed points
of $S$, the horospheres are permuted by $c$, while the elements
$a,b$ leaves each of them invariant.

Each $H_s$ is diffeomorphic to ${\bf R}^3$, and moreover
foliated by planes each of which is invariant by
the subgroup generated by $a,b$ in $S$, which is
isomorphic to ${\Bbb Z}^2$.
Indeed a leaf, $L_u$, of the foliation is the product of
a horosphere, which is diffeomorphic to ${\bf R}$,
in the first ${\bf H}^2$ in $X$
and a horosphere in the second ${\bf H}^2$.
By the definition of $\sigma$, $L_u$ is invariant by $\langle a,b \rangle $.
It follows that for each $s$, $H_s/S$ is
homeomorphic to $M$, and
$N=X/S$ is homeomorphic to $M \times {\bf R}$.
$N$ is foliated by $H_s/S, s \in {\bf R}$, which is
a product.

A simply connected, complete Riemannian
manifold whose sectional curvature is
non-positive is called an {\it Hadamard manifold}.
Hadamard manifolds are proper CAT(0) spaces.
It is easy to see that $S$ does not
act properly on any Hadamard manifold, $X$, of dimension $3$
by isometries.
Indeed if it did, then the action would be free
and the quotient $X/S$ is a manifold. It has to be
closed since otherwise, the cohomological
dimension of $S$ would be less than $3$, the dimension of
$X$, which is impossible. But if the action is co-compact,
then the action has to be by semi-simple isometries,
which is impossible because the elements $a,b$ are
forced to be parabolic from the presentation of $S$.

We record the discussion.
\begin{prop}
Let
$$S=\langle a,b,c\,|\,ab=ba,cac^{-1}=a^2b,cbc^{-1}=ab \rangle .$$
$S$ is the fundamental group of a closed
$3$-manifold, $M$, which is a torus
bundle over a circle.
Then

\begin{enumerate}
\item
$\geomdim (S)=3$.
\item
$S$ acts properly, hence freely, by isometries
on $X={\bf H}^2 \times {\bf H}^2$
such that $X/S$ is homeomorphic
to $M \times {\bf R}$.
\item
$S$ does not act properly on any Hadamard manifold
of dimension $3$
by isometries.

\end{enumerate}
\end{prop}

\begin{question}\label{bundle}
Does $S$ act properly on some $3$-dimensional CAT(0) space
by isometries?
\end{question}

\begin{remark}\label{S}
If $S$ acts properly on a proper CAT(0) space,
$X$, by isometries, then by Theorem \ref{AB},
there exists a point $v \in X(\infty)$ with
$S(v)=v$ because $S$ is solvable, hence amenable.
Look at the action of $S$ on the set, ${\mathcal H}$, of
horospheres centered at $v$. Since the action
factors through the abelianization of $S$, both $a$ and $b$ act
trivially. Suppose that the dimension of $X$ is $3$. Then
$c$ has to act non-trivially on ${\mathcal H}$
since otherwise $S$ would act trivially on ${\mathcal H}$,
so that we can apply Prop \ref{cdG} and
get ${\rm cd}S \le 2$, a contradiction.
\end{remark}

A CAT(0) space, $X$, is called {\it visible} if
for any two distinct points, $p,q$, in the ideal boundary,
there exists a geodesic, $\gamma$, in $X$
with $\gamma(\infty)=p,\gamma(-\infty)=q$.
We know a partial answer to Question \ref{bundle}.

\begin{prop}\label{visible}
$S$ does not act properly on any $3$-dimensional proper CAT(0)
space which is visible.
\end{prop}

\begin{proof}
Suppose such an action did exist on $X$ which is
visible. Then, as we said in Remark \ref{S},
there is a point, $v \in X(\infty)$, with $S(v)=v$.
Remark that since $X$ is visible,
a parabolic isometry has a unique fixed point in $X(\infty)$.
This is a standard fact (cf.\ \cite{FNS}).

We first show that $c$ can not be parabolic. Indeed, if it was parabolic,
a fixed point of $c$ given by Prop \ref{parabolicfix}
is the point $v$, so that $c$ has to
leave each horosphere centered at $v$ invariant.
Therefore, looking at the action of $S$ on the set of
the horospheres,
we get a contradiction as in Remark \ref{S}.
Thus, $c$ is hyperbolic. Let $\gamma$
be an axis of $c$.
Then, either $\gamma(\infty)$ or $\gamma(-\infty)$ is $v$, since
otherwise $\gamma$ would bound a flat half plane in $X$
because $c(v)=c$, which is impossible since $X$ is visible.

In the following argument, without loss
of generality, we assume that the action
is a right action, and denote the result of the action
by a group element $g$ of a point $y$ by $yg$,
instead of $gy$ or $g(y)$. For example, the point
$yg$ is mapped by $h$ to $ygh$.

Suppose $\gamma(\infty)=v$.
There is a constant
$R$ such that the displacement function of $a$ on $X$, $d_a$,
is at most $R$ on $\gamma(t), t \ge0$,
because $a(\gamma(\infty))=\gamma(\infty)$.

Set $x=\gamma(0)$. Let $K$ be the set of elements $g \in B$ such
that $d(x,xg) \le R$.
Then, $d(x,xcac^{-1}) = d(xc,xca) \le R$, because $xc$ is on $\gamma(t),
t \ge 0$.
Therefore $cac^{-1} \in K$.
By the same reason, for all $n \ge 1$,
$c^n a c^{-n} \in K$. On the other hand, they are all
different elements. One can see it
by rewriting each of them as a (positive) word in $\langle a,b \rangle$
using the relators in the presentation of $S$.

Suppose $\gamma(-\infty)=v$.
Fix an integer $N >0$. Let $x_N=xc^{-N}$.
For each $0 \le n \le N$,
$d(x_N, x_N c^n a c^{-n}) \le R$,
because $x_N c^n$ is on $\gamma(t), t \ge 0$.
Therefore, there are at least $N+1$
group elements, $g$, with $d(x_N,x_N g) \le R$.
However, the number of those elements for each $x_N$
does not depend on $N$, and also finite, because
the points $x_N$ are in one orbit of
the $S$-action, which contains $x$.
Since $N$ is arbitrary, we get a contradiction.
\end{proof}

By the same argument, we can show the following
proposition. In the proof, the roles of the elements $c$ and $a$ of $S$
in the
proof of Prop \ref{visible} are replaced by
the elements $a$ and $b$ of $BS(1,m)$, respectively. We omit details.
This is a partial
answer to Question \ref{bs.question}
\begin{prop}\label{prop.bs.vis}
$BS(1,m), m \ge 2$,
does not act on a proper, visible CAT(0) space
of dimension $2$.
\end{prop}

\subsection{Products}
If one looks for a group with the gap
between the cohomological dimension and
a CAT(0) dimension
to be $2$, or bigger, natural candidates to
be considered are the product of groups
with the gap $1$, which we already found.
To prove a theorem in this regard
using induction, one may want to use
a splitting theorem. Among several versions
of those(e.g., 6.21 in \cite{BH}),
we quote the following one, which seems least restrictive.

\begin{thm}[Monod, \cite{Mo}]\label{Mo}
Suppose $X$ is a proper CAT(0) space.
Suppose $G=G_1 \times \cdots \times G_n$
is acting on $X$ by isometries.

Then, either there is a point in $X(\infty)$ fixed by $G$,
or there is a non-empty closed convex $G$-invariant subset
$Z \subset X$ which splits $G$-equivariantly, isometrically
as $Z=Z_1 \times \cdots \times Z_n$ with an isometric
action of $G_i$ on each $Z_i$.
\end{thm}

Although we do not know how to show a desired result, let's
see what we can say. For example, let $G$ be either the
group $A$ in Theorem \ref{BC2} or the group $B$ in Theorem \ref{bridson}.
Set $H=G\times G$. Then $\cd H=4$.
Suppose $H$ acts properly on a proper CAT(0) space, $X$, by
isometries.
Since $\cd H=4$, $\dim X \ge 4$.
We would like to show that $\dim X \ge 6$.

It is easy to see that $\dim X \not =4$.
Suppose $\dim X=4$.
Recall that $G$ contains a subgroup, $K$,
isomorphic to ${\Bbb Z}^2$.
Consider the action of $K \times K$ on $X$.
Since $K \times K \simeq {\Bbb Z}^4$, by
Proposition \ref{ss}, the action is by semi-simple
isometries.
Then, by the flat torus theorem (Theorem \ref{flattorus}),
there is a convex subspace in $X$ which is isometric to ${\bf E}^2$
and invariant by $K \times \{e\}$.
It is a standard fact (cf.\ 7.1, II in {BrH}) that
the union of those flats in $X$ is isometric to ${\bf E}^2 \times Y$
such that $Y$ is a proper CAT(0) space.
Since $\dim X =4$, $\dim Y \le 2$.
The group $G$ is acting properly on $Y$ by isometries,
because $\{ e \} \times G$ centralizes $K \times \{ e \}$.
But this is impossible by Theorems \ref{BC2}, \ref{bridson}.

Therefore, the critical case is when $\dim X=5$.
Suppose $\dim X=5$. Then we can show that $H$ has to fix a point
in $X(\infty)$, because, otherwise, by Theorem \ref{Mo}, $H$ has to
act on some product space $Z=Z_1 \times Z_2$ whose
dimension is at most $5$. Therefore either $Z_1$ or
$Z_2$ has dimension less than $2$, but it is impossible
since $G$ acts properly on it by isometries.
In conclusion, to show that $\dim X \ge 6$, we are left
with the case that $H$ has a common fixed point in $X(\infty)$.
Probably we also need to analyze
a group action with a fixed point in the ideal boundary
to deal with Questions \ref{bs.question}, \ref{bundle}.

\Addresses\recd

\end{document}

%% file: agtout.tex
\def\ifplaintex{\expandafter\ifx\csname documentclass\endcsname\relax}

\def\gtp{{\mathsurround=0pt\it $\cal G\mskip-2mu$eometry \&\ 
$\cal T\!\!$opology $\cal P\!$ublications}}  

\def\recd{{\small Received:\qua\receiveddate\ifx\reviseddate\relax
\else\qquad Revised:\qua\reviseddate\fi\par}} 

\def\lognumber#1{\def\thelognumber{#1}}
\def\volumenumber#1{\def\thevolumenumber{#1}}
\def\volumeyear#1{\def\thevolumeyear{#1}}
\def\papernumber#1{\def\thepapernumber{#1}}
\def\pagenumbers#1#2{\def\startpage{#1}\def\finishpage{#2}}
\def\published#1{\def\publishdate{#1}}

\def\received#1{\def\receiveddate{#1}}
\def\revised#1{\def\reviseddate{#1}}
\def\accepted#1{\def\accepteddate{#1}}

\def\asciiemail#1{\def\theasciiemail{#1}}

\long\def\asciiabstract#1{\long\def\theasciiabstract{#1}}

\let\\\par\let\thelognumber\relax\let\thevolumenumber\relax
\let\thepapernumber\relax\let\thevolumeyear\relax\let\startpage\relax
\let\finishpage\relax\let\publishdate\relax\let\receiveddate\relax
\let\reviseddate\relax\let\accepteddate\relax\let\theasciititle\relax
\let\theasciiauthors\relax
\let\theasciiabstract\relax

\let\theasciiemail\relax

\ifplaintex
\font\logobig=cmssbx10 scaled 3836
\font\logomed=cmssbx10 scaled 2557
\else
\font\logobig=cmssbx10 scaled 4200
\font\logomed=cmssbx10 scaled 2800
\fi

\long\def\makeagttitle{   
\count0=\startpage
\agt\hfill      
\hbox to 45truept{\vbox to 0pt{\vglue -13truept{\logomed A\kern -.37em{\logobig 
T}\kern -.38em G}\vss}\hss}
\break
{\small Volume \thevolumenumber\ (\thevolumeyear)
\startpage--\finishpage\nl
Published: \publishdate}

\vglue .25truein

{\parskip=0pt\leftskip 0pt plus
1fil\def\\{\par\smallskip}{\Large\bf\thetitle}\par\medskip} \vglue
0.05truein

{\parskip=0pt\leftskip 0pt plus 1fil\def\\{\par}{\sc\theauthors}
\par\medskip}%
 
\vglue 0.03truein 

{\small\leftskip 25truept\rightskip 25truept{\bf Abstract}\stdspace\theabstract

{\bf AMS Classification}\stdspace\theprimaryclass
\ifx\thesecondaryclass\relax\else; \thesecondaryclass\fi\par
{\bf Keywords}\stdspace \thekeywords\par}\vglue 7truept

}   

\ifplaintex
\font\phead=cmsl9 scaled 950
\font\pnum=cmbx10 scaled 913
\font\pfoot=cmsl9 scaled 950
\headline{\vbox to 0pt{\vskip -4.5mm\line{\small\phead\ifnum
\count0=\startpage ISSN 1472-2739 (on-line) 1472-2747 (printed)
\hfill {\pnum\folio}\else\ifodd\count0\def\\{ }%
\ifx\theshorttitle\relax\thetitle\else\theshorttitle\fi\hfill{\pnum\folio}
\else\def\\{ and }{\pnum\folio}\hfill\ifx\theshortauthors\relax\theauthors
\else\theshortauthors\fi\fi\fi}\vss}}
\footline{\vbox to 0pt{\vglue 0mm\line{\small\pfoot\ifnum\count0=\startpage
\copyright\ \gtp\hfill\else
\agt, Volume \thevolumenumber\ (\thevolumeyear)\hfill\fi}\vss}}
\else
\font=cmsl9 scaled 1050
\font\lnum=cmbx10 
\font=cmsl9 scaled 1050
\makeatletter
\def\@oddhead{{\small\lhead\ifnum\count0=\startpage ISSN 1472-2739 
(on-line) 1472-2747 (printed)\hfill {\lnum\number\count0}\else\ifodd\count0
\def\\{ }\ifx\theshorttitle\relax \thetitle \else\theshorttitle\fi\hfill
{\lnum\number\count0}\else\def\\{ and }{\lnum\number\count0}
\hfill\ifx\theshortauthors\relax 
\theauthors\else\theshortauthors\fi\fi\fi}}\def\@evenhead{\@oddhead}
\def\@oddfoot{\small\lfoot\ifnum\count0=\startpage\copyright\ \gtp\hfill\else
\agt, Volume \thevolumenumber\ (\thevolumeyear)\hfill\fi}
\def\@evenfoot{\@oddfoot}
\makeatother
\fi
\let\maketitlepage\makeagttitle
\let\makeshorttitle\maketitlepage
\let\maketitle\maketitlepage

\newwrite\gtoutfile
\long\gdef\makeheadfile{  
{\def\\{, }\def\s{ }
\immediate\openout\gtoutfile head.xxx
\immediate\write\gtoutfile{Proxy-for: \ifx\theasciiauthors\relax
\theauthors\else\theasciiauthors\fi\s<\ifx\theasciiemail\relax\theemail\else\theasciiemail\fi>}
\immediate\write\gtoutfile{\noexpand\\}
\immediate\write\gtoutfile{Authors: \ifx\theasciiauthors\relax
\theauthors\else\theasciiauthors\fi}
{\def\\{ }\immediate\write\gtoutfile{Title: \ifx\theasciititle\relax
\thetitle\else\theasciititle\fi}}
\immediate\write\gtoutfile{Subj-class: GT or SG, GR etc}
\immediate\write\gtoutfile{MSC-class: \theprimaryclass\ifx\thesecondaryclass\relax\else, \thesecondaryclass\fi}
\immediate\write\gtoutfile{Journal-ref: Algebr. Geom. Topol. \thevolumenumber\s
(\thevolumeyear) \startpage-\finishpage}
\immediate\write\gtoutfile{Comments: Published by Algebraic and
Geometric Topology at}
\immediate\write\gtoutfile{\s\s\s  http://www.maths.warwick.ac.uk/agt/AGTVol\thevolumenumber/agt-\thevolumenumber-\thepapernumber.abs.html}
\immediate\write\gtoutfile{\noexpand\\}
\immediate\write\gtoutfile{}
\ifx\theasciiabstract\relax
\immediate\write\gtoutfile{\theabstract}\else
\immediate\write\gtoutfile{\theasciiabstract}\fi
\immediate\write\gtoutfile{}
\immediate\write\gtoutfile{\noexpand\\}
\immediate\write\gtoutfile{}
\immediate\closeout\gtoutfile}}  

\def\maketitlepage{\makeagttitle\makeheadfile}
\let\makeshorttitle\maketitlepage
\let\maketitle\maketitlepage

%% file: agt-4-38.bbl
\begin{thebibliography}{CheGr}
\itemsep 1.5pt plus 1pt

\bibitem[AB]{AB}
{\bf Scott Adams}, {\bf Werner Ballmann},
Amenable isometry groups of Hadamard spaces.
Math. Ann. 312 (1998), 183--195.
\MR{1645958}

\bibitem[B]{BB}
{\bf Werner Ballmann},
``Lectures on spaces of nonpositive curvature".
DMV Seminar, 25. Birkhauser, 1995
\MR{1377265}

\bibitem[BGS]{BGS}
{\bf Werner Ballmann}, {\bf Mikhael Gromov}, {\bf Victor Schroeder},
Manifolds of nonpositive curvature.
Progress in Mathematics, 61.
Birkhauser (1985).
\MR{0823981}

\bibitem[BeM]{BM}
{\bf Mladen Bestvina}, {\bf Geoffrey Mess},
The boundary of negatively curved groups.
J. Amer. Math. Soc. 4 (1991), no. 3, 469-481.
\MR{1096169}

\bibitem[BraC]{BC}
{\bf Noel Brady}, {\bf John Crisp},
Two-Dimensional Artin Groups with CAT(0) Dimension Three,
in ``Proceedings of the Conference on Geometric and Combinatorial Group
Theory, Part I (Haifa, 2000)'',
Geometriae Dedicata 94, No1 (2002), 185-214.
\MR{1950878}


\bibitem[Brad]{TBr}
{\bf Thomas Brady},
Complexes of nonpositive curvature for extensions of
$F\sb 2$ by $ Z$. Topology Appl. 63 (1995), no. 3, 267--275.
\MR{1334311}


\bibitem[Bri]{Bri}
{\bf Martin R. Bridson},
Length functions, curvature and the dimension of discrete groups.
Math. Res. Lett. 8 (2001), no. 4, 557-567.
\MR{1851271} 

\bibitem[BriH]{BH}
{\bf Martin R. Bridson}, {\bf Andre Haefliger},
``Metric spaces of non-positive curvature".
Grundlehren der Mathematischen Wissenschaften 319.
Springer, 1999.
\MR{1744486}

\bibitem[Br]{Br}
{\bf Kenneth S. Brown},
``Cohomology of groups".
Graduate Texts in Mathematics, 87.
Springer, 1982.
\MR{0672956} 

\bibitem[BuBuI]{BuBuI}
{\bf Dmitri Burago}, {\bf Yuri Burago}, {\bf Sergei Ivanov},
A course in metric geometry. Graduate Studies in Mathematics,
33. AMS, 2001.
\MR{1835418} 

\bibitem[ChaD]{CD}
{\bf Ruth Charney}, {\bf Michael W. Davis},
Finite $K(\pi, 1)$s for Artin groups.
in ``Prospects in topology", 110-124,
Ann. of Math. Stud. 138,
Princeton Univ. Press, 1995.
\MR{1368655} 

\bibitem[CheGr]{CG}
{\bf S.S.Chen}, {\bf L.Greenberg},
Hyperbolic spaces, in ``Contributions to Analysis",
Academic Press, 49-87, (1974).
\MR{0377765}

\bibitem[DeV]{DV}
{\bf Satya Deo}, {\bf K. Varadarajan},
Discrete groups and discontinuous actions.
Rocky Mountain J. Math. 27 (1997), no. 2, 559--583.
\MR{1466157} 

\bibitem[FNS]{FNS}
{\bf Koji Fujiwara}, {\bf Koichi Nagano}, {\bf Takashi Shioya},
Fixed point sets of parabolic isometries of CAT($0$)-spaces.
preprint. 2004.

\bibitem[G]{G}
{\bf M.Gromov},
Asymptotic invariants of infinite groups.
in ``Geometric group theory, Vol. 2 (Sussex 1991)",
1--295, LMS Lecture Note Ser. 182,
Cambridge Univ. Press, Cambridge, (1993).
\MR{1253544} 

\bibitem[Kapo]{Kapo}
{\bf Michael Kapovich},
Hyperbolic manifolds and discrete groups.
Progress in Mathematics, 183.
Birkhauser, 2001.
\MR{1792613}

\bibitem[Mo]{Mo}
{\bf N.Monod},
Superrigidity for irreducible lattices and geometric splitting.
preprint, Dec 2003.

\bibitem[Na]{Ng:dim}
{\bf J. Nagata}, ``Modern dimension theory", revised edition,
Sigma Series in Pure Mathematics, 2.
Heldermann Verlag, 1983.
\MR{0715431} 

\bibitem[P]{P}
{\bf A.R.Pears},
Dimension theory of general spaces.
Cambridge University Press, 1975.
\MR{0394604} 

\bibitem[Di]{Di}
{\bf Tammo. tom Dieck},
Transformation groups. de Gruyter Studies in Mathematics, 8.
Walter de Gruyter, 1987.
\MR{0889050} 

\end{thebibliography}
